\documentclass[a4paper,12pt]{article}
\usepackage{graphicx}
\usepackage{a4wide}
\usepackage[english]{babel}
\usepackage{amsfonts}
\usepackage{amsmath}
\usepackage{amssymb}
\usepackage{dsfont}

\newtheorem{theorem}{Theorem}

\allowdisplaybreaks

\date{}

\newcommand {\E} {\mathbb{E}}

\newcommand {\al} {\alpha}

\newcommand {\ve} {\varepsilon}

\numberwithin{equation}{section} \numberwithin{theorem}{section}
\numberwithin{lemma}{section} 
\numberwithin{corollary}{section}
\makeatletter\def\blfootnote{\xdef\@thefnmark{}\@footnotetext}\makeatother

\def\varep{\varepsilon}

\begin{document}
\title{\bf On permutation-invariance of limit theorems}
\author{I.\ Berkes\footnote{ Graz University of Technology,
Institute of Statistics, Kopernikusgasse 24, 8010 Graz, Austria.
\mbox{e-mail}: \texttt{berkes@tugraz.at}. Research supported by
FWF grants P24302-N18, W1230 and OTKA grant K 106814.} \, and R.\
Tichy\footnote{Graz University of Technology, Institute of
Mathematics A, Steyrergasse 30, 8010 Graz, Austria. \mbox{e-mail}:
\texttt{tichy@tugraz.at}. Research supported by FWF grants
P2304-N18, W1230 and SFB project F5510.}}

\maketitle \vskip0.5cm

\abstract{By a classical principle of probability theory,
sufficiently thin subsequences of general sequences of random
variables behave like i.i.d.\ sequences. This observation not only
explains the remarkable properties of lacunary trigonometric series,
but also provides a powerful tool in many areas of analysis, such
the theory of orthogonal series and Banach space theory. In contrast
to i.i.d.\ sequences, however, the probabilistic structure of
lacunary sequences is not permutation-invariant and the analytic
properties of such sequences can change after
rearrangement. In a previous paper we showed that
permutation-invariance of subsequences of the trigonometric system
and related function systems is connected with Diophantine
properties of the index sequence. In this paper we will study
permutation-invariance of subsequences of general r.v.\ sequences.}



\newpage
\section{Introduction}

It is known that sufficiently thin subsequences of general r.v.\
sequences behave like i.i.d.\  sequences. For example, R\'ev\'esz
\cite{re} showed that if a sequence $(X_n)$ of r.v.'s  satisfies
$\sup_n EX_n^2<\infty$, then one can find a subsequence $(X_{n_k})$
and a r.v.\ $X\in L^2$ such that $\sum_{k=1}^\infty c_k (X_{n_k}-X)$
converges a.s.\ provided $\sum_{k=1}^\infty c_k^2<\infty$. Under the
same condition, Gaposhkin \cite{gap1966}, \cite{gap1972} and
Chatterji \cite{chaclt}, \cite{chalil} proved that there exists a
subsequence $(X_{n_k})$ and r.v.'s  $X\in L^2$, $Y\in L^1$, $Y \ge
0$ such that
\begin{equation}\label{CLTm}
\frac{1}{\sqrt N} \sum_{k\le N} (X_{n_k} - X)
\overset{d}{\longrightarrow} N(0,Y)
\end{equation}
and
\begin{equation}\label{LILm}
\limsup_{N\to\infty} \, \frac{1}{\sqrt{2N \log\log N}} \sum_{k\le N}
(X_{n_k} - X) = Y^{1/2} \qquad \textup{a.s.} .
\end{equation}
Here $N(0, Y)$ denotes the distribution of the r.v.\ $Y^{1/2} \zeta$,
where $\zeta$ is a standard normal  r.v.\ independent of $Y$.
Koml\'os \cite{ko} showed that if $\sup_n E|X_n|<\infty$, then there exists  a
subsequence $(X_{n_k})$ and a r.v. $X\in L^1$ such that
$$\lim_{N\to\infty} \frac{1}{N} \sum_{k=1}^N X_{n_k}=X \qquad \text{a.s.} .$$
Chatterji \cite{chalp} showed that if $\sup_n E|X_n|^p<\infty$ where $0<p<2$, then
 the conclusion of the previous theorem can be changed to
$$
\lim_{N\to\infty} \frac{1}{N^{1/p}} \sum_{k=1}^N (X_{n_k}-X)=0
\qquad \text{a.s.}
$$
for some $X\in L^p$. Note the randomization in all these examples:
the  role of the mean and variance of the subsequence $(X_{n_k})$ is
played by random variables $X$, $Y$. For further limit theorems for
subsequences of general r.v.\ sequences and for the history of the
topic until 1966, see Gaposhkin \cite{gap1966}.

Since the asymptotic properties of an i.i.d.~sequence do not change
if we permute its terms, it is natural to expect that limit theorems
for lacunary subsequences of general r.v.~sequences remain valid
after any permutation of their terms.
This is, however, not the case. By classical results of Salem and
Zygmund \cite{sz1947}, \cite{sz1950} and Erd\H{o}s and G\'al
\cite{eg}, under the Hadamard gap condition
\begin{equation}\label{q}
n_{k+1}/n_k \ge q >1 \qquad k=1, 2, \ldots
\end{equation}
the sequence $(\sin 2\pi n_kx)$ satisfies
\begin{equation}\label{CLT}
\frac{1}{\sqrt{N/2}}\sum_{k=1}^N \sin 2\pi
n_kx\overset{d}{\longrightarrow}N(0, 1)
\end{equation}
and
\begin{equation}\label{LIL}
\limsup_{N\to\infty} \frac{1}{\sqrt{N\log\log N}}\sum_{k=1}^N \sin
2\pi n_kx=1 \qquad \text{a.s.}
\end{equation}
with respect to the probability space $( (0, 1), {\cal B}, \mu)$,
where $\mu$ denotes the Lebesgue measure. Erd\H{o}s \cite{er62} and
Takahashi \cite{tak72} proved that (\ref{CLT}), (\ref{LIL}) remain
valid under the weaker gap condition
\begin{equation} \label{ergap}
n_{k+1}/n_k \ge 1+ck^{-\alpha},  \qquad k=1, 2, \ldots
\end{equation}
for $0<\alpha<1/2$ and that for $\alpha=1/2$ this becomes false. As it was shown
in \cite{abt2}, \cite{abt3}, under the Hadamard gap condition (\ref{q})
the CLT (\ref{CLT}) and the LIL (\ref{LIL}) are permutation-invariant,
i.e.\ they remain valid after any permutation of the sequence
$(n_k)$, but this generally fails under the gap condition
(\ref{ergap}).
Similar results hold
for lacunary sequences $f(n_kx)$, where $f$ is a measurable function
satisfying
\begin{equation}\label{f}
f(x+1)=f(x), \qquad \int_0^1 f(x)\, dx=0, \qquad \int_0^1 f^2(x)\,
dx<\infty.
\end{equation} In this case, assuming the Hadamard gap
condition (\ref{q}), the validity of the CLT
\begin{equation}\label{CLTf}
\frac{1}{\sqrt{N}}\sum_{k=1}^N
f(n_kx)\overset{d}{\longrightarrow}N(0, \sigma^2)
\end{equation}
and of its permuted version depend on the number of solutions of the
Diophantine equation
\begin{equation}\label{dio}
a n_k + bn_\ell=c, \qquad 1\le k, \ell \le N.
\end{equation}
 As shown in  \cite{abt1}, \cite{abt2}, \cite{abt3}, a sharp condition for the CLT is that the number of
solutions of (\ref{dio}) is $o(N)$ for any fixed nonzero $a, b, c$, while
the permuted CLT requires the stronger bound $O(1)$ for the number
of solutions.

Permutation-invariance of limit theorems becomes a particularly difficult problem
for parametric limit theorems, e.g.\ for limit theorems containing arbitrary
coefficients. By a classical result of Menshov \cite{men}, from
every orthonormal system $(f_n)$ one can select a subsequence $(f_{n_k})$ which
is a convergence system, i.e.\ the series $\sum_{k=1}^\infty c_k f_{n_k}$ converges
almost everywhere provided $\sum_{k=1}^\infty c_k^2<\infty$. The question of whether
a subsequence $(f_{n_k})$ exists such that this property remains valid after any permutation
of $(f_{n_k})$ (i.e., by the standard terminology,  $(f_{n_k})$ is an unconditional
convergence system) remained open for nearly 40 years until it was answered in the
affirmative by Koml\'os \cite{ko2}. For another proof see Aldous \cite{ald}. The problem
of whether every orthonormal system can be rearranged  to become a convergence system is
still open; for a partial result see Garsia \cite{gar}.  Kolmogorov showed (see \cite{km})
that there exists an $f\in L^2(0,1)$ whose Fourier series, suitably permuted, diverges a.e.
But even though the Rademacher-Menshov convergence theorem yields a sharp a.e.\ convergence
criterion for orthonormal series, there is no similar complete result for rearranged
trigonometric series.

The previous results show that permutation-invariance of limit
theorems lies substantially deeper than that of the
original theorems and raise the question of which limit theorems hold in
a permutation-invariant form for lacunary sequences. In this paper we will
prove the surprising fact that, in a sense to be made precise, {\it all}
nonparametric distributional limit theorems for i.i.d.\ random variables
hold for lacunary subsequences $(f_{n_k})$ of general r.v.\ sequences in
a permutation-invariant form provided that the subsequence is sufficiently
thin, i.e.\ the gaps of the sequence (depending on the limit theorem) grow
sufficiently rapidly. We will deduce this result
from a general structure theorem for lacunary sequences proved in \cite{bp}
stating that sufficiently thin subsequences of any tight sequence of random
variables are nearly exchangeable. While this idea is simple and elementary,
formulating our results is somewhat technical and requires some preparations
in Section 2. The proof of our theorem will be given in Section 3.

\section{Main result}

We start with a formal definition of the concept "weak limit theorem".
Let $\cal M$
denote the set of all probability measures on $\mathbb R$ and $\varrho$ the
Prohorov metric on $\cal M$ defined by
\begin{align*}\varrho (\nu ,\lambda ) &= \inf \bigl\{ \varep > 0 : \nu (A)
\le \lambda (A^\varep )+\varep \ \hbox{ and} \cr
&\qquad \lambda (A) \le \nu (A^\varep) +\varep \ \hbox{ for all Borel
sets } \ A\subset \mathbb{R} \bigr\}.
\end{align*}
Here
$$A^\varep=\{x \in \mathbb{R}: |x-y|<\varep \ \text{for some} \ y\in A\}$$
denotes the open $\varep$-neighborhood of $A$.
A random measure is a measurable map from a probability space to
$\cal M$. The following definition is due to Aldous \cite{ald}.

\bigskip\noindent
{\bf Definition}. {\sl A weak limit theorem of i.i.d.~random
variables is a system
$$T =(f_1, f_2, \ldots, S, \{ G_\mu, \mu\in S\})$$
where

\smallskip\noindent
(a) $S$ is a Borel subset of $\cal M$;

\smallskip\noindent
(b) For each $k \ge 1$, $f_k = f_k(x_1, x_2, \ldots, \mu)$ is a
continuous function on ${\mathbb R}^\infty \times \cal M$,
satisfying the Lipschitz condition
$$|f_k(x_1, x_2, \ldots, \mu) - f_k (x'_1, x'_2, \ldots,
\mu)| \le \sum_{i = 1}^\infty c_{k, i} |x_i - x'_i|$$
where $0\le c_{k, i}\le 1$ and $\lim_{k\to\infty} c_{k, i}=0$ for all $i$;

\smallskip\noindent
(c) For each $\mu\in S$, $G_\mu$ is a probability
distribution on $\mathbb{R}$  such that the function $\mu \to G_\mu$
is  measurable (with respect to the Borel $\sigma$-fields in $S$ and
$\cal M$);

\medskip\noindent
and

\medskip\noindent
(d) If $\mu\in S$ and $X_1, X_2,\ldots$ are
independent r.v.'s with common distribution $\mu$ then
\begin{equation}\label{(3)}
f_k(X_1, X_2, \ldots, \mu)
\overset{d}{\longrightarrow} G_\mu \quad \text{as } k\to\infty.
\end{equation}
}

\bigskip
For example, the central limit theorem corresponds to
$$
S=\{\mu \in {\cal M}: \int x^2 d\mu (x)<\infty \}, \qquad G_\mu= N(0, \text{Var}\, \mu),
$$
$$
f_k (x_1, x_2, \ldots, \mu)= (x_1 + \ldots + x_k - k \cdot E\mu) /
\sqrt k, \qquad  c_{k, i}=k^{-1/2}I_{ \{i\le k\}}.
$$
The theorem itself is expressed by (\ref{(3)}).

\medskip
Using the terminology of \cite{bero}, we call a sequence $(X_n)$ of random variables
{\it determining}
if it has a limit distribution relative to any set $A$ in the
probability space with $P(A) > 0$, i.e.~for any $A\subset\Omega$
with $P(A) > 0$ there exists a distribution function $F_A$ such
that
$$\lim\limits_{n \to\infty} P(X_n < t\mid A) = F_A(t)$$
for all continuity points $t$ of $F_A$. Here $P(\cdot |A)$ denotes conditional probability
given $A$. (This concept is the same as that of stable convergence, introduced by
R\'enyi \cite{renyi}; our terminology follows that of functional analysis.) By an
extension of the Helly-Bray theorem (see \cite{bero}), every tight sequence of r.v.'s contains
a determining subsequence.
As is shown in \cite{ald}, \cite{bero}, for any determining sequence $(X_n)$
there exists a random measure $\tilde{\mu}$ (i.e.\ a measurable
map from the underlying probability space $(\Omega,\cal F, \cal
P)$ to $\cal M$) such that for any $A$ with $P(A) > 0$ and any
continuity point $t$ of $F_A$ we have
\begin{equation}\label{(4)} F_A(t) =
\E_A(\tilde{\mu}(-\infty, t)) \end{equation} where $\E_A$ denotes
conditional expectation given $A$. We call $\tilde{\mu}$ the {\it
limit random measure\/} of $(X_n)$.

The following result is Aldous' celebrated subsequence theorem \cite{ald}.

\bigskip\noindent
\begin{theorem} \label{permsubald} Let $(X_n)$ be a determining sequence with
limit random measure $\tilde{\mu}$. Let $T = (f_1, f_2, \ldots,
S, \{ G_\mu, \mu\in S\} )$ be a weak limit theorem and
assume $P(\tilde{\mu} \in S) = 1$. Then there exists a
subsequence $(X_{n_k})$ such that
\begin{equation} \label{pald}
f_k(X_{n_1}, X_{n_2}, \ldots, \tilde{\mu} ) \overset{d}{\longrightarrow} \int G_{\tilde{\mu}} dP
.
\end{equation}
\end{theorem}

\bigskip
In case of the CLT formalized above, assuming $\sup_n \E
X_n^2<+\infty$ implies easily that $\tilde{\mu}$ has finite variance
almost surely and thus denoting its mean and variance by $X$ and
$Y$, respectively, we see that the integral in (\ref{pald}) is the
distribution $N(0, Y)$. Hence (\ref{pald}) states in the present case
that
$$\frac{1}{\sqrt{N}} \sum_{k=1}^N (X_{n_k}-X) \overset{d}{\longrightarrow} N(0, Y)$$
which is exactly the CLT of Chatterji \cite{chaclt} and  Gaposhkin
\cite{gap1972} formulated in the Introduction. Theorem \ref{permsubald}
shows that a similar subsequence theorem holds for any weak limit theorem of
i.i.d.\ random variables. For a version of this result for strong (a.s.)
limit theorems, we refer to Aldous \cite{ald}.


In what follows, we change the technical conditions on $f_k$
in the definition of weak limit theorems slightly, leading to a class more
convenient for our purposes.

\bigskip\noindent
{\bf Definition.} {\sl The limit theorem $T =(f_1, f_2, \ldots, S,
 \{ G_\mu, \mu\in S\})$ is called {\it regular} if there exist two sequences
 $p_k\le q_k$ of positive integers tending to $+\infty$ and a sequence $\omega_k\to +\infty$ such that

\smallskip\noindent
 (i)  $f_k(x_1, x_2, \ldots, \mu)$ depends only on $x_{p_k}, \ldots, x_{q_k}, \mu$

\smallskip\noindent
 (ii) $f_k$ satisfies the Lipschitz condition
 \begin{equation} \label{lip}
|f_k(x_{p_k}, \ldots, x_{q_k}, \mu) - f_k (x'_{p_k},\ldots,
x'_{q_k}, \mu')| \le \frac{1}{\omega_k} \sum^{q_k}_{i = p_k} |x_i - x'_i|^\alpha
+ \varrho^* (\mu,\mu')
\end{equation}
for some $0 < \alpha \le 1$ where $\varrho^*$ is a metric on $S$
generating  the same topology as the Prohorov metric $\varrho$.}

\medskip
Thus in this case the function $f_k$ depends only on a finite
segment $x_{p_k}, \ldots x_{q_k}$  of the variables  $x_1, x_2,\ldots$ .
On the role of $\varrho^*$ see \cite{ald}. The above
definition brings out clearly the crucial feature of
limit theorems, namely the fact that the validity of the theorem does not
depend on finitely many terms of $(X_n)$, while the original
definition assumes only that the dependence of $f_k(X_1, X_2,
\ldots)$ on any fixed variable $X_j$ of the sequence is weak if $k$
is large. However, there is very little difference between these
assumptions. For example, the central limit theorem can be
formalized by either of the functions
$$
f_k (x_1, \ldots, x_k, \mu)= (x_1 + \ldots + x_k - k \cdot E\mu) /
\sqrt k
$$
and
$$
f^*_k(x_{[k^{1/4}]}, \ldots, x_k, \mu) =(x_{[k^{1/4} ]} +\ldots +
x_k - k \cdot E\mu)/\sqrt k
$$
of which the second leads to a regular limit theorem with the Wasserstein
metric
$$\varrho^* (\mu,\mu')=\left(\int_0^1 |F_\mu^{-1} (x) -F_{\mu'}^{-1}(x)|^2dx\right)^{1/2},$$
where $F_\mu, F_{\mu'}$ denote the distribution function of $\mu$ and $\mu'$,
respectively. Under bounded second moments, the contribution of the first
$k^{1/4}$ terms
in the normed sum defining $f_k$ are irrelevant and thus  we can always
switch from $f_k$ to $f^*_k$ and back again. The same procedure applies in
the general case.

We are now in a position to formulate the main result of our paper.

\bigskip\noindent
\begin{theorem} \label{permsub} Let $(X_n)$ be a determining sequence
with limit random measure $\tilde{\mu}$. Let $T = (f_1, f_2,
\ldots, S, \{ G_\mu, \mu\in S\} )$ be a regular weak limit theorem
and assume that $P(\tilde{\mu} \in S) = 1$. Then there exists a
subsequence $(X_{n_k}) = (Y_k)$ such that for any permutation
$(Y'_k) $ of $(Y_k)$ we have
\begin{equation} \label{p}
f_k(Y'_1, Y'_2, \ldots, \tilde{\mu} ) \overset{d}{\longrightarrow} \int G_{\tilde{\mu}}
dP .
\end{equation}
\end{theorem}

Note that we assumed the regularity of the limit theorem, but as we
pointed  out before, this is no restriction of generality.

The limit theorem $T$ in Theorem \ref{permsub} is non-parametric, i.e.\
the function $f_k$ depends on $x_1, x_2, \ldots$ and $\mu$, but on no additional parameters.
A simple example of a parametric distributional limit theorem is the weighted CLT, where
$$
f_k= A_k^{-1} \sum_{j=1}^k a_j (x_j-E\mu), \qquad A_k=\left(\sum_{j=1}^k a_j^2\right)^{1/2}.
$$
For any fixed coefficient sequence $(a_k)$ this defines a nonparametric limit theorem $T$ and Theorem
\ref{permsub} applies, but the selected subsequence $(X_{n_k})$ depends on $(a_k)$. As the
discussion above shows, in the case of a parametric limit theorem $T$ deciding whether a universal
subsequence $(X_{n_k})$ working for all parameters is generally a very difficult problem; an
example of a limit theorem where such a choice is impossible is given in \cite{gura}. For this reason,
in the present paper we deal only with nonparametric limit theorems.

In Aldous \cite{ald} a formalization of strong limit theorems is also given and the analogue of Theorem
\ref{permsubald} is proved. Using a reformulation of strong limit theorems as a sequence of probability
inequalities as given in \cite{be1985}, \cite{bp}, a version of our Theorem \ref{permsub} can be given
for a subclass of limit theorems considered in \cite{ald}. We also mention that for a more limited
class of weak limit theorems Theorem \ref{permsub} was proved in \cite{pe}.

\section{Proof of Theorem \ref{permsub}.}

To simplify the
formulas, let $f_k(\mu)$ denote, for any $\mu\in S$, the
distribution of the random variable $f_k(\xi_1, \xi_2, \ldots, \mu)$
where $\xi_1, \xi_2, \ldots$ are independent r.v.'s with common
distribution $\mu$. The following statements are easy to verify:
\medskip
\item{(A)} If $\varrho(\mu,\nu) \le \ve$ then $\varrho(f_k(\mu),
f_k(\nu)) \le \ve^\al q_k + \varrho^*(\mu,\nu)$ where $\al, q_k$ and
$\varrho^*$ are the quantities appearing in (\ref{lip}).

\item{(B)} Let $\mu_1, \ldots, \mu_r$ and $\nu_1,\ldots\nu_r$ be
probability distributions,  further let $c_1,\ldots, c_r$ be
nonnegative numbers with $\sum_{i=1}^r c_i = 1$. Assume that the sum of
those $c_i$'s such that $\varrho(\mu_i,\nu_i) \ge \ve$ is at most
$\ve$. Then the Prohorov distance between $\sum\limits^r_{i = 1} c_i
\mu_i$ and $\sum\limits^r_{i = 1} c_i \nu_i$ is at most $2\ve$.

\item{(C)} Let $\tilde{\mu}$ and $\tilde{\nu}$ be random measures
(i.e. measurable maps from a probability space $(\Omega,\cal F, P)$
to $\cal M$) such that $P(\varrho(\tilde{\mu}, \tilde{\nu}) \ge \ve)
\le \ve$. Then the Prohorov distance between $\int \tilde{\mu} dP$ and
$\int \tilde{\nu} dP$ is $\le 2\ve$.

\medskip
To prove statement (A) note that if $\varrho(\mu,\nu) \le \ve$ then
by a theorem of Strassen \cite{str}  there exist, on some
probability space, r.v.'s $\xi$ and $\eta$ with distribution $\mu$
and $\nu$ such that $P(|\xi - \eta| \ge \ve) \le \ve$. On a
larger probability space, let $(\xi_n, \eta_n)$ $(n = 1,2,\ldots)$ be
independent random vectors distributed as $(\xi,\eta)$. Clearly
$P(|\xi_i - \eta_i| \ge \ve) \le \ve$ $(i = 1,2,\ldots)$ and thus
using (\ref{lip}) we see that
$$\big|f_k (\xi_{p_k},\ldots, \xi_{q_k}, \mu) -
f_k(\eta_{p_k},\ldots, \eta_{q_k},\nu) \big| \le \ve^\al q_k +
\varrho^*(\mu,\nu)$$ except on a set with probability $\le \ve q_k
\le \ve^\al q_k$, proving (A). (Clearly we can assume $0<\ve\le 1$ and
that in the definition of regular limit theorems we have $\omega_k\ge 1$ for all $k$.)
Statements (B) and (C) are almost evident, (B) is a special case of (C).

To prove our theorem, let $(X_n)$ be a determining sequence of
r.v.'s with limit random measure $\tilde{\mu}$.  Then $(X_n)$ is
tight, i.e. $\sup_j P(|X_j| \ge t) \to 0$ as $t\to\infty$. Since
$\omega_k\to +\infty$, we can choose a nondecreasing sequence $(r_k)$ of
integers tending to $+\infty$ so slowly that
\begin{equation}\label{(2)}
r_k \le \min(p_k - 1, \omega^{1/4}_k)
\end{equation}
and
\begin{equation}\label{(13)} \sup_j P \Big(|X_j| \ge \frac{1}{2}
\omega^{1/(4\al)}_k\Big) \le \frac{1}{2} r_k^{-2} \qquad (k\ge
1).
\end{equation} Let $(\ve_k)$ tend to $0$
monotonically and so rapidly that
\begin{equation}\label{14}
\ve^\al_{r_k} q_k \le k^{-1} .
\end{equation} Using the structure
theorem \cite[Theorem 2]{bp}, it follows that there exists a subsequence
$(X_{n_k})$ and a sequence $(X'_k)$ of r.v.'s such that
\begin{equation}\label{(15)} |X_{n_k} - X'_k| = O(2^{-k}) \quad
\textup{a.s.}
\end{equation}
and $X'_k$ has the following properties:\\

\noindent \ \ (A$_1$) \ \ Each $X'_k$ takes only finitely many
values

\smallskip\noindent
\ \ (B$_1$)\ \ $\sigma\{ X'_1\} \subset \sigma \{ X'_2\}
\subset\ldots$

\smallskip\noindent
\ \ (C$_1$)\ \ For each $k\ge 1$ the atoms of the finite
$\sigma$-field $\sigma \{ X'_{r_k} \}$ can be divided into two
classes $\Gamma_1$ and $\Gamma_2$ so that
\begin{equation}\label{(16)}
\sum_{A\in\Gamma_1} P(A) \le \ve_{r_k}
\end{equation}
and for any $A\in \Gamma_2$ there exist i.i.d.r.v.'s $\{ Z^{(A)}_j,
j = r_k + 1, r_k + 2, \ldots \}$ defined on $A$ with distribution
function $F_A$ such that
\begin{equation}\label{(17)}
P_A \big(| X'_j - Z^{(A)}_j | \ge \ve_{r_k} \big) \le \ve_{r_k}
\quad j = r_k + 1, r_k + 2,\ldots.\\
\end{equation}

\noindent Here $F_A$ denotes the limit distribution of $(X_n)$ on
the set $A$ (which exists since $(X_n)$ is determining) and $P_A$
denotes conditional probability with respect to $A$.
\medskip
Let $\tilde{\mu}_n$ denote the random measure defined by
$\tilde{\mu}_n (B) = \E (\tilde{\mu} (B) \mid X_n')$. By Lemma~7 of
\cite{bp}  we have $\tilde{\mu}_n \overset{d}{\longrightarrow} \tilde{\mu}$ a.s. and thus
by passing to a further subsequence of $(X_{n_k})$ we can also
assume that
\begin{equation}\label{(18)}
P\{ \varrho(\tilde{\mu}_n, \tilde{\mu}) \ge \ve_n\} \le \ve_n
\end{equation}
\begin{equation}\label{(19)}
P\{ \varrho^*(\tilde{\mu}_n, \tilde{\mu}) \ge \ve_n\} \le \ve_n .
\end{equation}
We show that the last obtained subsequence $(X_{n_k})$ satisfies the
conclusion of the theorem. In view of (\ref{lip}) and (\ref{(15)}),  $X_{n_k}$
and $X'_k$ are interchangeable in the statement of the theorem and
thus it suffices to prove that if $(X'_k)$ satisfies statements
(A$_1$), (B$_1$), (C$_1$) above then for any permutation $(Y'_k)$ of
$(X'_k)$ we have (\ref{p}). To verify this, note that by  (\ref{lip})
and (\ref{(17)}) we have
\begin{equation}\label{(20)}
\gathered P_A \big\{ \big| f_k(X'_{i_1},\ldots, X'_{i_\ell}, \mu_A)
- f_k(Z^{(A)}_{i_1}, \ldots, Z^{(A)}_{i_\ell}, \mu_A)\big| \ge
\ve^\al_{r_k} q_k \big\}\\
\hskip3cm \le \ve^\al_{r_k} q_k \quad \qquad A\in\Gamma_2
\endgathered
\end{equation}
where $\ell = q_k - p_k + 1$, $i_1,\ldots, i_\ell$ are different
integers $>r_k$ and $\mu_A$ is the probability measure corresponding
to $F_A$. (Note that we do not assume here  $i_1<\ldots<i_\ell$~;
the vectors $(X'_{i_1} \ldots,  X'_{i_\ell})$ and $(Z^{(A)}_{i_1}, \ldots, Z^{(A)}_{i_\ell})$ are close to each other coordinatewise, i.e.\ for any order of $i_1, \ldots, i_\ell$. Since the $Z_j^{(A)}$ are i.i.d., the distribution of the vector $(Z^{(A)}_{i_1}, \ldots, Z^{(A)}_{i_\ell})$ is permutation-invariant, providing an explanation for the phenomenon described in Theorem \ref{permsub}.)
Since (\ref{(20)}) is valid for all $A\in\Gamma_2$ and
$\mu_A$ in (\ref{(20)}) is identical to $\tilde{\mu}_{r_k}$ on $A$
(see Lemma 6 of \cite{be1985}), using (\ref{(16)}), (\ref{(20)}) and
statement (B) at the beginning of the proof we get
\begin{equation}\label{(21)}
\varrho\big( f_k(X'_{i_1},\ldots, X'_{i_\ell}, \tilde{\mu}_{r_k}),
\sum_A f_k(\mu_A) P(A)\big) \le 2\ve^\al_{r_k} q_k
\end{equation}
where the sum is extended for all atoms $A$ of $\sigma\{ X'_{r_k}\}$
and a r.v. in a Prohorov distance is meant as its distribution. Next
we show that (\ref{(21)}) remains valid, with the right hand side
increased by $r^{-1}_k$, if $i_1,\ldots, i_\ell$, $\ell = q_k - p_k
+ 1$, are arbitrary different positive integers (not necessarily
$>r_k$). Indeed, remove from $X'_{i_1},\ldots, X'_{i_\ell}$ those
whose index is $\le r_k$ and replace them with (different) $X'_j$'s
with $j > \max(r_k, i_1,\ldots, i_\ell)$. This means that we change
$f_k(X'_{i_1},\ldots, X'_{i_\ell}, \tilde{\mu}_{r_k})$ at most at
$r_k$ locations and at each such position we replace an $X'_\mu$ by
an $X'_\nu$ where $\mu \le r_k$ and $\nu > r_k$. By (\ref{lip}), $f_k$
changes at most by
$$\frac{1}{\omega_k} \sum |X'_\mu - X'_\nu|^\al =: W$$
where the sum has $\le r_k$ terms. Using (\ref{(2)}), (\ref{(13)})
we get
\begin{align*}
& P(|W| \ge r^{-1}_k) \le P(|W| \ge \omega^{-1/2}_k) =
P\left(\sum|X'_\mu - X'_\nu|^\al \ge \omega^{1/2}_k\right)\\
&\le \sum P\Big( |X'_\mu - X'_\nu| \ge \Big(
\frac{\omega^{1/2}_k}{r_k} \Big)^{1/\al}\Big) \le 2r_k \cdot
\sup\nolimits_j P \Big( |X'_j| \ge \frac{1}{2}
\omega^{1/(4\al)}_k\Big)\le r^{-1}_k
\end{align*}
and thus the above changes increase the left hand side of
(\ref{(21)}) by at most $r^{-1}_k$, i.e.
\begin{equation}\label{(22)}
\varrho\Big( f_k (X'_{i_1},\ldots, X'_{i_\ell},\tilde{\mu}_{r_k}),
\sum_A f_k(\mu_A) P(A)\Big) \le 2 \ve^\al_{r_k} q_k + r^{-1}_k
\end{equation}
for any different positive integers $i_1,\ldots, i_\ell$, $\ell =
q_k - p_k + 1$. Changing $\tilde{\mu}_{r_k}$ into $\tilde{\mu}$ will
change $f_k (X'_{i_1},\ldots, X'_{i_\ell},\tilde{\mu}_{r_k})$ on the
left hand side of (\ref{(22)}) by at most
$\ve_{r_k}$, except on a set of probability $\le \ve_{r_k}$ (see
(\ref{(19)}) and (\ref{lip})) and thus the left hand side of (\ref{(22)})
changes by at most $\ve_{r_k}$. Thus observing that the sum $\sum_A
f_k(\mu_A) P(A)$ in (\ref{(22)}) equals $\int f_k(\tilde{\mu}_{r_k})
dP$, we proved
the following\\

\medskip\noindent
{\bf Proposition}. Let $(X^*_k)$ be any permutation of $(X'_k)$.
Then
$$\varrho\big( f_k(X^*_{p_k},\ldots, X^*_{q_k},\tilde{\mu}), \int
f_k(\tilde{\mu}_{r_k}) dP\big) \le 3 \ve^\al_{r_k} q_k + r^{-1}_k
.$$

\medskip
To complete the proof of our theorem it suffices to show that the
Prohorov distance of any two of the distributions
\begin{equation}\label{(23)}
\int f_k(\tilde{\mu}_{r_k}) dP \qquad \int f_k(\tilde{\mu})dP \qquad
\int G_{\tilde{\mu}} dP \end{equation} tends to zero as
$k\to\infty$. To verify this observe first that (\ref{(18)}),
(\ref{(19)}) and statement (A) at the beginning of  the proof imply
that the Prohorov distance of $f_k(\tilde{\mu}_{r_k})$ and
$f_k(\tilde{\mu})$ is $\le \ve^\al_{r_k} q_k + \ve_{r_k}$, except on
a set with probability $\le \ve^\al_{r_k} q_k + \ve_{r_k}$ and thus
by statement (C) and (\ref{(4)}) the Prohorov distance of the first
two distributions in (\ref{(23)}) is $\le 2(\ve^\al_{r_k} q_k +
\ve_{r_k}) \le 4k^{-1}$. On the other hand, the validity of
$f_k(\mu)\overset{d}\longrightarrow G_\mu$ for any $\mu\in S$ (which is a part of the
definition of a weak limit theorem) and $P(\tilde{\mu} \in S) = 1$
imply $\varrho(f_k(\tilde{\mu}), G_{\tilde{\mu}}) \to 0$ a.s. and
thus there exists a numerical sequence $\delta_k \downarrow 0$ such
that
$$P\{ \varrho(f_k(\tilde{\mu}), G_{\tilde{\mu}}) \ge \delta_k\} \le
\delta_k \qquad (k = 1,2,\ldots).$$ Thus by statement (C) above we
get that the Prohorov distance of the second and third distribution
in (\ref{(23)}) is $\le 2\delta_k$. This completes the proof
of Theorem \ref{permsub}.\\

\end{document}